\documentclass[10pt]{article}
\usepackage{amssymb,amsfonts}
\usepackage{amsmath}
\usepackage{amscd,amssymb,amsthm}
\usepackage{graphics}
\usepackage{graphicx}
\usepackage{hyperref}

\graphicspath{ {./images/} }
\hypersetup{
    colorlinks=true,citecolor=blue
  }
\newtheorem{theorem}{Theorem}[section]

\newtheorem{lemma}[theorem]{Lemma}

\newtheorem{remark}[theorem]{Remark}

\begin{document}

\title{Bifurcation diagrams in a class of Kolmogorov systems}
\author{G. Moza, C. Lazureanu\thanks{
Department of Mathematics, Politehnica University of Timisoara, Romania;
gheorghe.moza@upt.ro}, F. Munteanu, C. Sterbeti, A. Florea\thanks{%
Department of Applied Mathematics, University of Craiova, Romania} }
\date{}
\maketitle

\begin{abstract}
We study a two-dimensional Kolmogorov system when its two parameters vary in a small neighbourhood of the value $0.$ The local behavior of the system is described in terms of bifurcation diagrams.
\end{abstract}


\section{Introduction}

We consider in this work the following Kolmogorov class of systems

\begin{equation}
\begin{tabular}{lll}
$\frac{dx}{d\tau }$ & $=$ & $2x\left( \mu
_{1}+p_{11}x+p_{12}y+p_{13}xy+s_{1}x^{2}\right) $ \\ 
&  &  \\ 
$\frac{dy}{d\tau }$ & $=$ & $2y\left( \mu
_{2}+p_{21}x+p_{22}y+p_{23}x^{2}+s_{2}y^{2}\right) $%
\end{tabular}%
,  \label{s1}
\end{equation}%
where $\mu =\left( \mu _{1},\mu _{2}\right) \in 
\mathbb{R}
^{2},$ $p_{ij}=p_{ij}\left( \mu \right) $ and $s_{i}=s_{i}\left( \mu \right) 
$ are smooth functions of their arguments, $i=1,2,$ $j=1,2,3,$ $\left\vert
\mu \right\vert =\sqrt{\mu _{1}^{2}+\mu _{2}^{2}}$ is sufficiently small, $%
0\leq \left\vert \mu \right\vert \ll 1,$ and $p_{12}\left( 0\right)
p_{22}\left( 0\right) \neq 0.$ More difficult is the case $p_{12}\left(
0\right) p_{22}\left( 0\right) =0,$ which will be tackled in another work.

In general, differential systems of the form 
\begin{equation}
\dot{x}_{i}=x_{i}\cdot P_{i}\left( x_{1},...,x_{n}\right) ,  \label{ko1}
\end{equation}%
$i=1,...,n,$ $x_{i}=x_{i}\left( t\right) ,$ $P_{i}$ smooth functions, are
called Kolmogorov systems \cite{jl1}. Their applications can be found in
many areas, such as biology \cite{bel}, \cite{rafi}, fluid dynamics
(turbulence) \cite{bul} and plasma physics \cite{laval}. For example, a
Kolmogorov three-dimensional system has been recently studied in \cite{adak}
as a model for parasites dynamics. If $P_{i}$ are polynomials, then (\ref%
{ko1}) is a Lotka--Volterra system which is widely used in modeling
interacting species of predator-prey type arising, for example, in biology 
\cite{bra}, \cite{free} and ecology \cite{kot}, \cite{may}, \cite{yang}.

Since many of these applications use positive variables, we will study the
system (\ref{s1}) only in the first quadrant $Q1,$ that is, when $x\geq 0$
and $y\geq 0.$ However, a similar study can be performed for the other
quadrants. It is worth mentioning that the axes $x=0$ and $y=0$ are
invariant manifolds with respect to system's flow, thus, an orbit starting
in the first quadrant cannot cross any of the two axes and travel in other
quadrants. It remains forever (for all $t$ where it is defined) in the first
quadrant. This aspect is important in applications because, if a model of
type (\ref{ko1}) has practical relevance only when its variables $x_{i}$ are
positive (and this is the case in most of the models of this type), then, an
orbit starting in the zone of practical relevance will never leave this zone
to enter a zone which does not present interest, characterized by $%
x_{i}\left( t\right) <0$ for some $i\in \left\{ 1,..,n\right\} $ and $t\in 
\mathbb{R}
.$ In other words, if the starting state of a model lies in a zone of
practical relevance, then the evolution of that state remains in the same
zone for any time.

The paper is organized as it follows. After the introduction, in section 2
we present the main analytical results which describe the behavior of the
system. In section 3 we present the bifurcation diagrams arising in the
system. An in-depth study of the system is performed in this section, by
finding the curves which separate a node from a focus. Such curves are in
general ignored in typical bifurcation analysis but they are highly
appreciated in practical models used in biology and engineering, for
example. Differentiating between a node and a focus, twenty six different
phase portraits are obtained in the bifurcation diagrams. Conclusive remarks
end the article.

\section{Analytical properties of the system}

Depending on the signs of $p_{12}\left( 0\right) $ and $p_{22}\left(
0\right),$ the system leads to two different forms, one corresponding to 1) $%
p_{12}\left( 0\right) p_{22}\left( 0\right) <0$ and the other to 2) $%
p_{12}\left( 0\right) p_{22}\left( 0\right) >0.$

Indeed, let us consider first 1) with $p_{12}\left( 0\right) <0$ and $%
p_{22}\left( 0\right) >0.$ Define the changes 
\begin{equation}
\xi _{1}=-p_{12}\left( \mu \right) x,\xi _{2}=p_{22}\left( \mu \right) y%
\text{ and }t=2\tau .  \label{ec7}
\end{equation}%
The transformation of variables $\left( x,y\right) \longmapsto \left( \xi
_{1},\xi _{2}\right) $ is well defined for all $\left\vert \mu \right\vert $
small enough and is nonsingular because $p_{12}\left( 0\right) p_{22}\left(
0\right) \neq 0.$

It follows from (\ref{ec7}) that the system (\ref{s1}) is locally
topologically equivalent near the origin to 
\begin{equation}
\left\{ 
\begin{tabular}{lll}
$\frac{d\xi _{1}}{dt}$ & $=$ & $\xi _{1}\left( \mu _{1}-\theta \left( \mu
\right) \xi _{1}+\gamma \left( \mu \right) \xi _{2}-M\left( \mu \right) \xi
_{1}\xi _{2}+N\left( \mu \right) \xi _{1}^{2}\right) $ \\ 
&  &  \\ 
$\frac{d\xi _{2}}{dt}$ & $=$ & $\xi _{2}\left( \mu _{2}-\delta \left( \mu
\right) \xi _{1}+\xi _{2}+S\left( \mu \right) \xi _{1}^{2}+P\left( \mu
\right) \xi _{2}^{2}\right) $%
\end{tabular}%
\right. ,  \label{s2}
\end{equation}%
where the coefficients $\theta \left( \mu \right) =\frac{p_{11}\left( \mu
\right) }{p_{12}\left( \mu \right) },$ $\gamma \left( \mu \right) =\frac{%
p_{12}\left( \mu \right) }{p_{22}\left( \mu \right) },$ $M\left( \mu \right)
=\frac{p_{13}\left( \mu \right) }{p_{12}\left( \mu \right) p_{22}\left( \mu
\right) },$ $N\left( \mu \right) =\frac{s_{1}\left( \mu \right) }{%
p_{12}^{2}\left( \mu \right) },$ $\delta \left( \mu \right) =\frac{%
p_{21}\left( \mu \right) }{p_{12}\left( \mu \right) },$ $S\left( \mu \right)
=\frac{p_{23}\left( \mu \right) }{p_{12}^{2}\left( \mu \right) }$ and $%
P\left( \mu \right) =\frac{s_{2}\left( \mu \right) }{p_{22}^{2}\left( \mu
\right) }$ are smooth functions in $\mu .$ In what follows, these
coefficients are needed only at $\mu =0.$ In order to save symbols, we
denote by $\theta \left( 0\right) =\theta ,$ $\gamma \left( 0\right) =\gamma
,$ $\delta \left( 0\right) =\delta ,$ $N=N\left( 0\right) $ and so on.
Notice that $\theta =\frac{p_{11}\left( 0\right) }{p_{12}\left( 0\right) },$ 
$\gamma =\frac{p_{12}\left( 0\right) }{p_{22}\left( 0\right) }$ and $\delta =%
\frac{p_{21}\left( 0\right) }{p_{12}\left( 0\right) }$ are well-defined from 
$p_{12}\left( 0\right) p_{22}\left( 0\right) \neq 0.$ Moreover, $\gamma <0$
while $\theta \delta $ can be $0.$

\begin{remark}
\label{rem1} If $p_{12}\left( 0\right) >0$ and $p_{22}\left( 0\right) <0,$
using the changes $\xi _{1}=p_{12}\left( \mu \right) x,$ $\xi
_{2}=-p_{22}\left( \mu \right) y$ and $t=-2\tau ,$ the system (\ref{s1})
reads 
\begin{equation*}
\left\{ 
\begin{tabular}{lll}
$\frac{d\xi _{1}}{dt}$ & $=$ & $\xi _{1}\left( -\mu _{1}-\theta \left( \mu
\right) \xi _{1}+\gamma \left( \mu \right) \xi _{2}+M\left( \mu \right) \xi
_{1}\xi _{2}-N\left( \mu \right) \xi _{1}^{2}\right) $ \\ 
&  &  \\ 
$\frac{d\xi _{2}}{dt}$ & $=$ & $\xi _{2}\left( -\mu _{2}-\delta \left( \mu
\right) \xi _{1}+\xi _{2}-S\left( \mu \right) \xi _{1}^{2}-P\left( \mu
\right) \xi _{2}^{2}\right) $%
\end{tabular}%
\right. ,
\end{equation*}%
which is of the same form as (\ref{s2}) if one denotes by $\widetilde{\mu }%
_{1,2}=-\mu _{1,2},$ $\widetilde{M}\left( \mu \right) =-M\left( \mu \right)
, $ $\widetilde{N}\left( \mu \right) =-N\left( \mu \right) ,$ $\widetilde{S}%
\left( \mu \right) =-S\left( \mu \right) $ and $\widetilde{P}\left( \mu
\right) =-P\left( \mu \right) .$
\end{remark}

In the case 2), assume first $p_{12}\left( 0\right) >0$ and $p_{22}\left(
0\right) >0.$ Then (\ref{s1}) is locally topologically equivalent near the
origin to 
\begin{equation}
\left\{ 
\begin{tabular}{lll}
$\frac{d\xi _{1}}{dt}$ & $=$ & $\xi _{1}\left( \mu _{1}+\theta \left( \mu
\right) \xi _{1}+\gamma \left( \mu \right) \xi _{2}+M\left( \mu \right) \xi
_{1}\xi _{2}+N\left( \mu \right) \xi _{1}^{2}\right) $ \\ 
&  &  \\ 
$\frac{d\xi _{2}}{dt}$ & $=$ & $\xi _{2}\left( \mu _{2}+\delta \left( \mu
\right) \xi _{1}+\xi _{2}+S\left( \mu \right) \xi _{1}^{2}+P\left( \mu
\right) \xi _{2}^{2}\right) $%
\end{tabular}%
\right. ,  \label{s3}
\end{equation}%
by using the changes $\xi _{1}=p_{12}\left( \mu \right) x,$ $\xi
_{2}=p_{22}\left( \mu \right) y$ and $t=2\tau .$ We notice also that, the
case $p_{12}\left( 0\right) <0$ and $p_{22}\left( 0\right) <0$ reduces to (%
\ref{s3}) by changes $\xi _{1}=-p_{12}\left( \mu \right) x,$ $\xi
_{2}=-p_{22}\left( \mu \right) y$ and $t=-2\tau ,$ using also new parameters
such as $\widetilde{\mu }_{1,2}=-\mu _{1,2}$ and renaming the coefficients $%
M,$ $N,$ $S$ and $P$ as needed.

\begin{remark}
\label{rem2} Since $\gamma \left( 0\right) >0$ in (\ref{s3}) and $\gamma
\left( 0\right) <0$ in (\ref{s2}), the two systems are not necessarily
locally topologically equivalent near the origin, thus, they should be
studied separately. In this work we tackle the first case, $p_{12}\left(
0\right) p_{22}\left( 0\right) <0.$
\end{remark}

\begin{remark}
\label{rem2b} The axes $\xi _{1}=0$ and $\xi _{2}=0$ in (\ref{s2}) are
invariant with respect to the system's flow and, by (\ref{ec7}), $\xi
_{1}\geq 0$ and $\xi _{2}\geq 0$ whenever $x\geq 0$ and $y\geq 0.$ It
follows that, the first quadrant of (\ref{s1}) is transformed by (\ref{ec7})
in the first quadrant of (\ref{s2}), denoted by $Q1.$ Thus, the new system (%
\ref{s2}) will be studied only in the first quadrant also.
\end{remark}

We approach in this work the first case, $p_{12}\left( 0\right) p_{22}\left(
0\right) <0,$ and more precisely the system (\ref{s2}) with $\theta \left(
0\right) \delta \left( 0\right) \neq 0$ and $\gamma \left( 0\right) <0.$

The first three equilibria of (\ref{s2}) are $O\left( 0,0\right) ,$ $%
E_{1}\left( \xi _{1},0\right) \text{ and }E_{2}\left( 0,\xi _{2}\right) ,$
where $\mu _{1}-\theta \xi _{1}+N\xi _{1}^{2}=0$ and $\mu _{2}+\xi _{2}+P\xi
_{2}^{2}=0,$ which, in their lowest terms become $\xi _{1}=\frac{1}{\theta }%
\mu _{1}+\allowbreak O\left( \mu _{1}^{2}\right) $ and $\xi _{2}=-\mu
_{2}+\allowbreak O\left( \mu _{2}^{2}\right) ,$ for all $\left\vert \mu
\right\vert $ sufficiently small.

The existence of a fourth equilibrium $E_{3}\left( \xi _{1},\xi _{2}\right) $
is ensured by the Implicit Function Theorem (IFT) applied to the algebraic
system 
\begin{equation}
\mu _{1}-\theta \xi _{1}+\gamma \xi _{2}-M\xi _{1}\xi _{2}+N\xi _{1}^{2}=0%
\text{ and }\mu _{2}-\delta \xi _{1}+\xi _{2}+S\xi _{1}^{2}+P\xi _{2}^{2}=0,
\label{as}
\end{equation}%
$\allowbreak $provided that 
\begin{equation}
\theta -\gamma \delta \neq 0.  \label{c1}
\end{equation}%
The coordinates of $E_{3}\left( \xi _{1},\xi _{2}\right) $ are of the form

\begin{equation}
\xi _{1}=\frac{\mu _{1}}{\theta -\gamma \delta }\left( 1+O\left( \left\vert
\mu \right\vert \right) \right) -\frac{\gamma \mu _{2}}{\theta -\gamma
\delta }\left( 1+O\left( \left\vert \mu \right\vert \right) \right) \text{
and }\xi _{2}=\frac{\delta \mu _{1}}{\theta -\gamma \delta }\left( 1+O\left(
\left\vert \mu \right\vert \right) \right) -\frac{\theta \mu _{2}}{\theta
-\gamma \delta }\left( 1+O\left( \left\vert \mu \right\vert \right) \right) .
\label{x12}
\end{equation}%
\qquad

We denote by $O\left( \left\vert \mu \right\vert ^{k}\right) =\sum_{i+j\geq
k}c_{ij}\mu _{1}^{i}\mu _{2}^{j}$ a Taylor series starting with terms of
order at least $k\geq 1.$

\begin{remark}
We are concerned in this work with local properties of the system (\ref{s2})
around the equilibrium points $O,$ $E_{1,2}$ and $E_{3}$ when $\left( \mu
_{1},\mu _{2}\right) $ moves in the parametric plane with $\left\vert \mu
\right\vert =\sqrt{\mu _{1}^{2}+\mu _{2}^{2}}$ sufficiently small.
\end{remark}

We state first the following result, which will be needed for studying the
existence of bifurcation curves. Denote by $D\subset 
\mathbb{R}
^{2}$ an open disc of center $\left( 0,0\right) $ and radius $\varepsilon >0$
sufficiently small, i.e. $D=\left\{ \left( x,y\right)
,x^{2}+y^{2}<\varepsilon \right\} $ is a small neighborhood of $\left(
0,0\right) .$

\begin{lemma}
\label{lemma2}Let $F,G:D\subset 
\mathbb{R}
^{2}\rightarrow 
\mathbb{R}
$ be two smooth functions of the form

a) $F\left( \mu _{1},\mu _{2}\right) =a\mu _{2}\left( 1+O\left( \left\vert
\mu \right\vert \right) \right) +b\mu _{1}\left( 1+O\left( \left\vert \mu
\right\vert \right) \right) ,$ $ab\neq 0,$ and

b) $G\left( \mu _{1},\mu _{2}\right) =a\mu _{2}^{2}\left( 1+O\left(
\left\vert \mu \right\vert \right) \right) -2b\mu _{1}\mu _{2}\left(
1+O\left( \left\vert \mu \right\vert \right) \right) +c\mu _{1}^{2}\left(
1+O\left( \left\vert \mu \right\vert \right) \right) ,$ $ac\neq 0.$ Then,
when $\left\vert \mu \right\vert $ is sufficiently small,

i) $F\left( \mu _{1},\mu _{2}\right) =0$ is a unique curve in the plane $\mu
_{1}\mu _{2},$ given by 
\begin{equation*}
\mu _{2}=-\frac{b}{a}\mu _{1}\left( 1+O\left( \mu _{1}\right) \right) ,
\end{equation*}

ii) if $\Delta =b^{2}-ac>0,$ then $G\left( \mu _{1},\mu _{2}\right) =0$ is a
union of two curves passing through the origin $O\left( 0,0\right) \ $and
having the forms 
\begin{equation}
\mu _{2}=e_{1}\mu _{1}\left( 1+O\left( \mu _{1}\right) \right) \text{ and }%
\mu _{2}=e_{2}\mu _{1}\left( 1+O\left( \mu _{1}\right) \right) ,
\label{c12a}
\end{equation}%
where $e_{1}=\frac{b+\sqrt{\Delta }}{a}$ and $e_{2}=\frac{b-\sqrt{\Delta }}{a%
},$ and

iii) if $\Delta =b^{2}-ac<0,$ then $G\left( \mu _{1},\mu _{2}\right) =0$ is
satisfied only at $\left( 0,0\right) $ and $sign\left( G\left( \mu _{1},\mu
_{2}\right) \right) =sign\left( a\right) $ for all $\left( \mu _{1},\mu
_{2}\right) \in D.$
\end{lemma}

PROOF. i) Since $F\left( 0,0\right) =0$ and $\frac{\partial F}{\partial \mu
_{2}}\left( 0,0\right) =a\neq 0,$ the IFT yields the conclusion. The
function $\mu _{2}=-\frac{b}{a}\mu _{1}\left( 1+O\left( \mu _{1}\right)
\right) $ exists and is unique for all $\left\vert \mu _{1}\right\vert $
sufficiently small.

ii) The proof is more involved in this case. $G\left( \mu _{1},\mu
_{2}\right) =0$ is equivalent to 
\begin{equation}
a\mu _{2}^{2}-2b\mu _{1}\mu _{2}\left( 1+O\left( \left\vert \mu \right\vert
\right) \right) +c\mu _{1}^{2}\left( 1+O\left( \left\vert \mu \right\vert
\right) \right) =0.  \label{gr2}
\end{equation}%
Notice that $\frac{1+O\left( \left\vert \mu \right\vert \right) }{1+O\left(
\left\vert \mu \right\vert \right) }=1+O\left( \left\vert \mu \right\vert
\right) .$ Then, solving for $\mu _{2}$ in (\ref{gr2}) when $\Delta >0,$ we
get an implicit equation 
\begin{equation*}
H_{\pm }\left( \mu _{1},\mu _{2}\right) \overset{def}{=}\mu _{2}-\frac{\mu
_{1}}{a}\left( b\left( 1+O\left( \left\vert \mu \right\vert \right) \right)
\pm \sqrt{\Delta }\left( 1+O\left( \left\vert \mu \right\vert \right)
\right) \right) =0.
\end{equation*}%
But $H_{\pm }\left( 0,0\right) =0$ and $\frac{\partial H_{\pm }}{\partial
\mu _{2}}\left( 0,0\right) =1\neq 0$ for all $O\left( \left\vert \mu
\right\vert \right) .$ Notice that $O\left( \left\vert \mu \right\vert
\right) =\sum_{i+j\geq 1}c_{ij}\mu _{1}^{i}\mu _{2}^{j}$ contains terms both
in $\mu _{1}$ but also in $\mu _{2}.$ The conclusion follows now from the
IFT and Taylor theorem. Notice that $e_{1,2}\neq 0$ from $ac\neq 0.$

iii) It is clear that, the implicit functions $H_{\pm }\left( \mu _{1},\mu
_{2}\right) $ do not exist when $\Delta <0.$ Thus $G\left( \mu _{1},\mu
_{2}\right) \neq 0$ for all $\left( \mu _{1},\mu _{2}\right) \neq \left(
0,0\right) $ with $\left\vert \mu \right\vert <\varepsilon .$ Moreover,
since $G\left( \mu _{1},\mu _{2}\right) $ is a quadratic equation in $\mu
_{2},$ it follows that $sign\left( G\left( \mu _{1},\mu _{2}\right) \right)
=sign\left[ a\left( 1+O\left( \left\vert \mu \right\vert \right) \right) %
\right] =sign\left( a\right) $ for $\left\vert \mu \right\vert $
sufficiently small. $\blacksquare $

\begin{remark}
If $b^{2}-ac=0,$ the discriminant of the equation $G\left( \mu _{1},\mu
_{2}\right) =0$ in $\mu _{2}$ is undefined, being of the form $\left\vert
\mu _{1}\right\vert \sqrt{O\left( \left\vert \mu \right\vert \right) }.$
Thus, the existence of curves of equation $G\left( \mu _{1},\mu _{2}\right)
=0$ cannot be decided. A further study could be performed in this case but
this is beyond the purpose of this article.
\end{remark}

The following bifurcation curves $T_{1}$ and $T_{2}$ are important for
obtaining bifurcation diagrams to describe the behavior of the system (\ref%
{s2}), namely 
\begin{equation}
T_{1}=\left\{ \left( \mu _{1},\mu _{2}\right) \in 
\mathbb{R}
^{2}\mid \mu _{2}=\frac{\delta }{\theta }\mu _{1}+O\left( \mu
_{1}^{2}\right) ,\frac{\mu _{1}}{\theta }>0\right\}  \label{T1}
\end{equation}%
and

\begin{equation}
T_{2}=\left\{ \left( \mu _{1},\mu _{2}\right) \in 
\mathbb{R}
^{2}\mid \mu _{2}=\frac{1}{\gamma }\mu _{1}+O\left( \mu _{1}^{2}\right) ,\mu
_{1}>0\right\} .  \label{T2}
\end{equation}%
$T_{1}$ is defined by $\xi _{2}=0$ while $T_{2}$ by $\xi _{1}=0,$ where $\xi
_{1,2}$ are given by (\ref{x12}) and correspond to $E_{3}\left( \xi _{1},\xi
_{2}\right) .$

\begin{remark}
\label{rem3}1) The eigenvalues of $E_{1}\left( \xi _{1},0\right) $ are $%
\lambda _{1}^{E_{1}}=\mu _{1}-2\theta \xi _{1}+3N\xi _{1}^{2}$ and $\lambda
_{2}^{E_{1}}=\mu _{2}-\delta \xi _{1}+S\xi _{1}^{2},$ which in their lowest
terms become $\lambda _{1}^{E_{1}}=-\mu _{1}+O\left( \mu _{1}^{2}\right) $
and $\lambda _{2}^{E_{1}}=\mu _{2}-\frac{1}{\theta }\delta \mu
_{1}+\allowbreak O\left( \mu _{1}^{2}\right) ,$ since $\xi _{1}=\frac{1}{%
\theta }\mu _{1}+\allowbreak O\left( \mu _{1}^{2}\right) .$

2) The eigenvalues of $E_{2}\left( 0,\xi _{2}\right) $ are $\lambda
_{1}^{E_{2}}=\mu _{1}+\gamma \xi _{2}$ and $\lambda _{2}^{E_{2}}=\mu
_{2}+2\xi _{2}+3P\xi _{2}^{2},$ that is, $\lambda _{1}^{E_{2}}=\mu
_{1}-\gamma \mu _{2}+O\left( \mu _{2}^{2}\right) $ and $\lambda
_{2}^{E_{2}}=-\mu _{2}+\allowbreak O\left( \mu _{2}^{2}\right) ,$ with $\xi
_{2}=-\mu _{2}+\allowbreak O\left( \mu _{2}^{2}\right) .$

3) The eigenvalues of $O$ are $\mu _{1,2}.$
\end{remark}

\begin{theorem}
\label{th0} Assume $\theta \delta \neq 0$ and $\theta -\gamma \delta \neq 0.$
Then $T_{1}$ is a transcritical bifurcation curve for all $\left\vert \mu
\right\vert $ sufficiently small.
\end{theorem}

\textit{Proof.} When $\left( \mu _{1},\mu _{2}\right) \in T_{1},$ the points 
$E_{3}\left( \xi _{1},\xi _{2}\right) $ and $E_{1}\left( \xi _{1},0\right) $
coalesce (collide), becoming a double non-hyperbolic singular point on $%
T_{1},$ with one eigenvalue $\lambda _{2}^{E_{1}}=0$ and the other $\lambda
_{1}^{E_{1}}=-\mu _{1}+O\left( \mu _{1}^{2}\right) .$ Indeed, $\xi _{2}=0$
in (\ref{as}) yields $\mu _{2}-\delta \xi _{1}+S\xi _{1}^{2}=0,$ thus, by
Remark \ref{rem3}, $\lambda _{2}^{E_{1}}=0.$ We notice that $E_{1}\in Q1$ if~%
$\frac{\mu _{1}}{\theta }>0,$ since $\xi _{1}=\frac{1}{\theta }\mu
_{1}+\allowbreak O\left( \mu _{1}^{2}\right) .$

We apply in the following Sotomayor's theorem (see \cite{perko}). To this
end, write first the system (\ref{s2}) in the form $\frac{d\xi }{dt}=f\left(
\xi ,\mu \right) ,$ where $\xi =\left( \xi _{1},\xi _{2}\right) ,$ $\mu
=\left( \mu _{1},\mu _{2}\right) $ and $f=\left( f_{1},f_{2}\right) .$

At $E_{1}\left( \xi _{1},0\right) $ with $\mu _{1}-\theta \xi _{1}+N\xi
_{1}^{2}=0,$ denote also by $\xi _{0}=\left( \xi _{1},0\right) $ and $\mu
_{0}=\left( \mu _{1},\mu _{2}\right) \in T_{1},$ with $\mu _{1}\neq 0.$
Since $E_{3}\left( \xi _{1},\xi _{2}\right) $ collides to $E_{1}\left( \xi
_{1},0\right) $ on $T_{1},$ by (\ref{as}), $\xi _{1}$ satisfies also $\mu
_{2}-\delta \xi _{1}+S\xi _{1}^{2}=0.$ Then, $f\left( \xi _{0},\mu
_{0}\right) =\left( 0,0\right) $ and the Jacobian matrix of $f$ at $\left(
\xi _{0},\mu _{0}\right) $ is $$A=Df\left( \xi _{0},\mu _{0}\right) =\left( 
\begin{array}{cc}
\theta \xi _{1}-2\mu _{1} & \xi _{1}\left( \gamma -M\xi _{1}\right) \\ 
0 & 0%
\end{array}%
\right) .$$ The eigenvector of $A,$ respectively $A^{T},$ corresponding to
the null eigenvalue $\lambda _{2}^{E_{1}}=0$ are $v=\left( 
\begin{array}{c}
\frac{\gamma \xi _{1}-M\xi _{1}^{2}}{2\mu _{1}-\theta \xi _{1}} \\ 
1%
\end{array}%
\right) ,$ respectively, $w=\allowbreak \left( 
\begin{array}{c}
0 \\ 
1%
\end{array}%
\right) .$ Notice that $v=\allowbreak \left( 
\begin{array}{c}
\frac{\gamma }{\theta }+O\left( \mu _{1}\right) \\ 
1%
\end{array}%
\right) $ is well-defined for all $\left\vert \mu \right\vert $ sufficiently
small with $\mu _{1}\neq 0.$ Assume $\mu _{1}$ is fixed with $\frac{\mu _{1}}{\theta }>0$ and consider $%
\mu _{2}$ as the bifurcation parameter. Denote by $f_{\mu _{2}}=\left( 
\begin{array}{c}
\frac{\partial f_{1}}{\partial \mu _{2}} \\ 
\frac{\partial f_{2}}{\partial \mu _{2}}%
\end{array}%
\right) $ and $D^{2}f\left( \xi ,\mu \right) \left( v,v\right) =\left( 
\begin{array}{c}
d^{2}f_{1}\left( \xi ,\mu \right) \left( v,v\right) \\ 
d^{2}f_{2}\left( \xi ,\mu \right) \left( v,v\right)%
\end{array}%
\right) ,$ where $d^{2}f_{1,2}\left( \xi ,\mu \right) \left( v,v\right) $
are the differentials of second order of the functions $f_{1,2}.$ Since $%
\frac{\partial f_{2}}{\partial \mu _{2}}\left( \xi _{0},\mu _{0}\right) =0,$
it follows that $C_{1}=w^{T}\cdot f_{\mu _{2}}\left( \xi _{0},\mu _{0}\right) =0.$
If $Df_{\mu _{2}}$ denotes the Jacobian matrix in variables $\xi _{1}$ and $%
\xi _{2}$ of $f_{\mu _{2}},$ we get 
\begin{equation*}
C_{2}=w^{T}\left[ Df_{\mu _{2}}\left( \xi _{0},\mu _{0}\right) \left(
v\right) \right] =\allowbreak \allowbreak 1+O\left( \mu _{1}\right) .
\end{equation*}
Finally, we find 
\begin{equation*}
C_{3}=w^{T}\left[ D^{2}f\left( \xi _{0},\mu _{0}\right) \left( v,v\right) %
\right] =2\frac{\theta -\gamma \delta }{\theta }+O\left( \mu _{1}\right) .
\end{equation*}
Thus, $C_{2}\neq 0$ and $C_{3}\neq 0$ for all $\left\vert \mu \right\vert $
sufficiently small, which complete the proof. $\blacksquare $

\begin{remark}
\label{rem3a}A similar scenario to $T_{1}$ occurs for $T_{2}.$ Indeed, if $%
\left( \mu _{1},\mu _{2}\right) \in T_{2},$ then $E_{3}\left( \xi _{1},\xi
_{2}\right) $ and $E_{2}\left( 0,\xi _{2}\right) $ coalesce. By (\ref{as}),
their joint eigenvalues are $\lambda _{1}^{E_{2}}=0$ and $\lambda
_{2}^{E_{2}}=-\mu _{2}+O\left( \mu _{2}^{2}\right) .$ $E_{2}\in Q1$ if $\mu
_{2}<0.$
\end{remark}

Denote by $X_{+}=\left\{ \left( \mu _{1},0\right) ,\mu _{1}>0\right\} ,$ $%
X_{-}=\left\{ \left( \mu _{1},0\right) ,\mu _{1}<0\right\} ,$ $Y_{+}=\left\{
\left( 0,\mu _{2}\right) ,\mu _{2}>0\right\} $ and $Y_{-}=\left\{ \left(
0,\mu _{2}\right) ,\mu _{2}<0\right\} ,$ the four half-axes in the
parametric plane $\mu _{1}\mu _{2}.$

On $X_{+}\cup X_{-},$ the point $O\left( 0,0\right) $ collides to $%
E_{2}\left( 0,\xi _{2}\right) $ and have the eigenvalues $0$ and $\mu _{1},$
while on $Y_{+}\cup Y_{-},$ $O\left( 0,0\right) $ collides to $E_{1}\left(
\xi _{1},0\right) $ and have the eigenvalues $0$ and $\mu _{2}.$

\begin{remark}
\label{rem3b}One can show similarly to $T_{1}$ that $T_{2},$ $X_{+},$ $%
X_{-}, $ $Y_{+}$ and $Y_{-}$ are all transcritical bifurcation curves.
Knowing that the system (\ref{s2}) undergoes a transcritical bifurcation
when the parameter $(\mu _{1},\mu _{2})$ crosses these curves is important
in determining the phase portraits on the curves, which represent borders of
different regions in the parametric plane. This will be exploited in the
next section.
\end{remark}

The characteristic polynomial at an equilibrium point $\left( \xi _{1},\xi
_{2}\right) $ of the system (\ref{s2}) has the form $P\left( \lambda \right)
=\lambda ^{2}-2p\allowbreak \lambda +L.$ When $\left( \xi _{1},\xi
_{2}\right) $ is an equilibrium point of type $E_{3},$ that is, its
coordinates satisfy (\ref{as}), then $p$ and $L$ at $\left( \xi _{1},\xi
_{2}\right) $ become%
\begin{equation}
p=\frac{1}{2}\xi _{2}-\frac{1}{2}\theta \xi _{1}+N\xi _{1}^{2}-\frac{1}{2}%
M\xi _{1}\xi _{2}+P\xi _{2}^{2},  \label{pg}
\end{equation}%
\begin{equation}
L=\allowbreak -\xi _{1}\xi _{2}\allowbreak \left( \allowbreak \theta -\gamma
\delta +c_{11}\xi _{1}+\allowbreak c_{12}\xi _{2}-4NP\xi _{1}\xi _{2}-2MS\xi
_{1}^{2}+2MP\xi _{2}^{2}\right) ,  \label{lg}
\end{equation}%
where $c_{11}=M\delta -2N+2S\gamma $ and $c_{12}=M+2P\theta .$

Consider that the eigenvalues of $E_{3}\left( \xi _{1},\xi _{2}\right) ,$
which are roots of $P\left( \lambda \right) ,$ are of the form 
\begin{equation}
\lambda _{1,2}=p\pm \sqrt{q}.  \label{lam1lam2}
\end{equation}%
Then $\lambda _{1}\lambda _{2}=L$ and $q=p^{2}-L\allowbreak $ lead to 
\begin{eqnarray}
q &=&\allowbreak \frac{1}{4}\left( \theta ^{2}\xi _{1}^{2}+2\xi _{1}\xi
_{2}\left( \theta -2\gamma \delta \right) +\xi _{2}^{2}\right) -N\theta \xi
_{1}^{3}+\allowbreak \frac{1}{2}\left( 2N+M\theta +2c_{11}\right)
\allowbreak \xi _{1}^{2}\xi _{2}+\frac{1}{2}\left( M+2P\theta \right) \xi
_{1}\xi _{2}^{2}  \notag \\
&&+P\xi _{2}^{3}+N^{2}\xi _{1}^{4}-\left( N+2S\right) M\xi _{1}^{3}\xi _{2}+%
\frac{1}{4}\left( M^{2}-8NP\right) \allowbreak \xi _{1}^{2}\xi
_{2}^{2}+MP\xi _{1}\xi _{2}^{3}+P^{2}\xi _{2}^{4}.  \label{qg}
\end{eqnarray}

Using the expressions of $\xi _{1,2}$ at $E_{3}$ given in (\ref{x12}), we
can write

\begin{equation*}
p=\allowbreak \frac{\delta -\theta }{2\left( \theta -\gamma \delta \right) }%
\mu _{1}\left( 1+O\left( \left\vert \mu \right\vert \right) \right) -\frac{%
\theta -\gamma \theta }{2\left( \theta -\gamma \delta \right) }\allowbreak
\mu _{2}\left( 1+O\left( \left\vert \mu \right\vert \right) \right) ,
\end{equation*}

\begin{equation}
q=\frac{1}{4\left( \theta -\gamma \delta \right) ^{2}}\left( a\mu
_{2}^{2}\left( 1+O\left( \left\vert \mu \right\vert \right) \right) -2b\mu
_{1}\allowbreak \mu _{2}\left( 1+O\left( \left\vert \mu \right\vert \right)
\right) +c\allowbreak \mu _{1}^{2}\left( 1+O\left( \left\vert \mu
\right\vert \right) \right) \right) ,  \label{q}
\end{equation}%
where 
\begin{equation}
a=\theta ^{2}\left( \gamma +1\right) ^{2}-4\theta \gamma ^{2}\delta
~,~b=-2\delta ^{2}\gamma ^{2}+\theta \left( \theta -\delta \right) \gamma
+\theta \left( \theta +\delta \right) \text{ and }c=\left( \theta +\delta
\right) ^{2}-4\gamma \delta ^{2}.  \label{coef}
\end{equation}%
Notice that $c\neq 0$ because $\theta \delta \neq 0$ and $\gamma <0.$ By (%
\ref{lg}), the product $\lambda _{1}\lambda _{2}$ is of the form 
\begin{equation}
\lambda _{1}\lambda _{2}=-\left( \theta -\gamma \delta \right) \xi _{1}\xi
_{2}\left( 1+O\left( \left\vert \mu \right\vert \right) \right) .
\label{l12}
\end{equation}

\begin{remark}
By Lemma \ref{lemma2}, we can approximate $E_{3},$ $p$ and $q$ by their
lowest terms, whenever they satisfy the Lemma's conditions. The qualitative
behavior of the system does not change by this approximation. For example,
the curve $\mu _{2}=e_{1}\mu _{1}\left( 1+O\left( \mu _{1}\right) \right) $
for $e_{1}\neq 0$ is approximated by $\mu _{2}=e_{1}\mu _{1}.$ In this
approximation, $sign\left[ G\left( \mu _{1},\mu _{2}\right) \right] =sign%
\left[ G_{0}\left( \mu _{1},\mu _{2}\right) \right] $ for $\left\vert \mu
\right\vert $ sufficiently small, where $G_{0}\left( \mu _{1},\mu
_{2}\right) =a\mu _{2}^{2}-2b\mu _{1}\mu _{2}+c\mu _{1}^{2}.$
\end{remark}

When $a\neq 0,$ the discriminant in its lowest terms of the equation $q=0$
in variable $\mu _{2}$ is 
\begin{equation}
\Delta =b^{2}-ac=\allowbreak -4\gamma \delta \left( \theta -\gamma \delta
\right) ^{3}\mu _{1}^{2}.  \label{delta}
\end{equation}

\begin{theorem}
\label{t1} 1) Assume $\theta -\gamma \delta <0.$ Then $E_{3}$ is stable
(node or focus) on $p<0$ and unstable on $p>0.$ In particular, if $\delta >0$
then $E_{3}$ is an unstable node.

2) Assume $\theta -\gamma \delta >0.$ Then $E_{3}$ is a saddle.
\end{theorem}

PROOF. It follows from (\ref{pg}) that $p>0$ if $\theta -\gamma \delta \neq
0 $ and $\theta <0,$ whenever $E_{3}\left( \xi _{1},\xi _{2}\right) $
satisfies $\xi _{1,2}>0.$ In its lowest terms, $p$ reads

\begin{equation*}
p=\allowbreak -\frac{1}{2\left( \theta -\gamma \delta \right) }\left( \theta
\left( \mu _{1}-\gamma \mu _{2}\right) -\left( \delta \mu _{1}-\theta \mu
_{2}\right) \right) .
\end{equation*}

1) When $\theta -\gamma \delta <0,$ $E_{3}\left( \frac{\mu _{1}-\gamma \mu
_{2}}{\theta -\gamma \delta },\frac{\delta \mu _{1}-\theta \mu _{2}}{\theta
-\gamma \delta }\right) $ exists in the first quadrant $Q1$ defined by $\xi
_{1}>0$ and $\xi _{2}>0,$\ whenever $\left( \mu _{1},\mu _{2}\right) $ lies
in the region 
\begin{equation*}
R_{1}=\left\{ \left( \mu _{1},\mu _{2}\right) \in 
\mathbb{R}
^{2}\mid \mu _{1}-\gamma \mu _{2}<0,\delta \mu _{1}-\theta \mu
_{2}<0\right\} .
\end{equation*}%
It follows from (\ref{l12}) that $\lambda _{1}\lambda _{2}>0$ on $R_{1}.$
Thus, $E_{3}$ is stable on $p<0$ (node or focus), respectively, unstable on $%
p>0$ (node or focus). Assume further:

i) $\delta >0.$ Then $\theta <\gamma \delta <0,$ $\Delta =\allowbreak
\allowbreak -4\gamma \delta \left( \theta -\gamma \delta \right) ^{3}<0,$ $%
a>0$ and $p>0,$ since $\gamma <0.$ Thus $q>0$ and $E_{3}$ is an unstable
node.

ii) $\delta <0.$ Then $\Delta >0$ and $q=0$ give rise to two bifurcation
curves, denoted by $C_{1}$ and $C_{2}$ and given in (\ref{c12a}), which
determine on $R_{1}$ when $E_{3}$ changes from a node to a focus or vice
versa; $\lambda _{1}\lambda _{2}>0.$

2) In the second case $\theta -\gamma \delta >0,$ $E_{3}\in Q1$ when $\left(
\mu _{1},\mu _{2}\right) $ lies in the region 
\begin{equation*}
R_{2}=\left\{ \left( \mu _{1},\mu _{2}\right) \in 
\mathbb{R}
^{2}\mid \mu _{1}-\gamma \mu _{2}>0,\delta \mu _{1}-\theta \mu
_{2}>0\right\} .
\end{equation*}%
From (\ref{l12}), we deduce $\lambda _{1}\lambda _{2}<0.$ Therefore $\lambda
_{1,2}\in \mathbb{R}$ and $q>0.$ We conclude that $E_{3}$ is a saddle. $%
\blacksquare $

\bigskip

Taking into account the behavior of $E_{3}$ described in the above theorem,
we can study easier the Hopf bifurcation in (\ref{s2}) at $E_{3}.$

\begin{theorem}
\label{t2} 1) If $\theta <0$ or $0<\gamma \delta <\theta $ or $\gamma \delta
<0<\theta ,$ the system (\ref{s2}) does not undergo a Hopf bifurcation at $%
E_{3}.$

2) A Hopf bifurcation occurs at $E_{3}$ provided that $0<\theta <\gamma
\delta .$
\end{theorem}

PROOF. 1) A Hopf bifurcation cannot occur at $E_{3}$ when $\theta <0$ either
on $R_{1}$ or $R_{2}$ because $p\neq 0.$ Assume further $\theta >0.$ Then $%
p=0$ along the curve 
\begin{equation}
H=\left\{ \left( \mu _{1},\mu _{2}\right) \left\vert \mu _{2}=\allowbreak
\allowbreak \mu _{1}\frac{\theta -\delta }{\theta \left( \gamma -1\right) }%
+O\left( \mu _{1}^{2}\right) ,\theta \neq \delta \right. \right\}  \label{h1}
\end{equation}%
and $q=\allowbreak \mu _{1}^{2}\frac{\theta -\gamma \delta }{\theta \left(
\gamma -1\right) ^{2}}+O\left( \mu _{1}^{3}\right) $ on $H;$ $\gamma <0.$ If
in addition $\gamma \delta <\theta ,$ then $q>0$ on $H.$ Thus, a Hopf
bifurcation cannot occur since the eigenvalues $\lambda _{1,2}$ are not
purely imaginary on $H.$

We still need to consider the case $\theta =\delta .$ Then $p$ in its lowest
terms is an expression of the form 
\begin{equation*}
p=-\frac{1}{2}\mu _{2}+\frac{1}{2}c_{n}\mu _{1}^{n}
\end{equation*}%
with $n\geq 2$ and $c_{n}\neq 0.$ Thus, $p=0$ leads to $\mu _{2}=c_{n}\mu
_{1}^{n}.$ Calculating now the expression of $q$ for $\mu _{2}=c_{n}\mu
_{1}^{n}$ we find $q=\allowbreak \frac{\mu _{1}^{2}}{1-\gamma }\left(
\allowbreak 1+O\left( \mu _{1}\right) \right) \allowbreak >0$ for $%
\left\vert \mu _{1}\right\vert $ small enough with $\mu _{1}\neq 0.$ It
concludes that a Hopf bifurcation can not occur at $\theta =\delta
\allowbreak .$

2) If $0<\theta <\gamma \delta $ then $\lambda _{1,2}=\pm \omega i$ on $H,$
with $\omega =\allowbreak \frac{\mu _{1}}{\theta \left( 1-\gamma \right) }%
\sqrt{\left( \gamma \delta -\theta \right) \theta }>0$ in their lowest
terms. One can show that $E_{3}$ is defined on $\mu _{1}>0$ in this case,
Figure \ref{d4} (VII). Since $\left. \frac{\partial p}{\partial \mu _{2}}%
\right\vert _{H}=\frac{\theta \left( \gamma -1\right) }{2\left( \theta
-\gamma \delta \right) \delta }\neq 0,$ a Hopf bifurcation (degenerate or
not) occurs at $E_{3}.$ We notice that $\theta \neq \delta $ in this case
because $\theta >0$ and $\delta <0.$

When the bifurcation is degenerate with respect to the first Lyapunov
coefficient (that is, $l_{1}\left( 0\right) =0$), other codim $k$
bifurcations may arise with $k\geq 2.$ In case of non-degeneracy of the Hopf
bifurcation, a periodic orbit (limit cycle) arises around $E_{3},$ but when
the bifurcation becomes degenerate (e.g. Bautin bifurcation) more limit
cycles are typically born around $E_{3}$ with the increasing of the
codimension. $\blacksquare $

\section{Bifurcation diagrams}

We perform in this section a complete analysis of the bifurcation diagrams
when the parameters $\delta ,\gamma ,\theta $ vary with $\theta \delta \neq
0,$ $\theta -\gamma \delta \neq 0$ and $\gamma <0.$ In order to obtain the
diagrams and distinct generic phase portraits of the system (\ref{s2}), we
will use the cases and results obtained in Theorems \ref{t1}-\ref{t2}.

\textbf{Case 1. }Consider the first inequality $\theta -\gamma \delta <0.$
More different cases arise here in order to obtain all configurations of
bifurcation diagrams.

\begin{figure}[h!]
\centering
\includegraphics[width=0.9\textwidth]{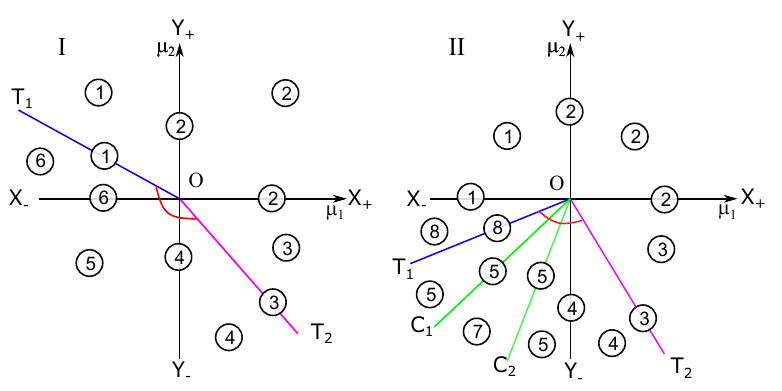}
\caption{Bifurcation diagrams of system (\protect\ref{s2}) for: $\protect%
\theta-\protect\gamma\protect\delta<0,$ $\protect\delta >0$ (I),
respectively, $\protect\theta-\protect\gamma\protect\delta<0,$ $\protect%
\delta <0,$ $\protect\theta<0,$ $a>0$ and $-1<\protect\gamma<0$ (II). }
\label{d1}
\end{figure}

\textbf{Case 1.a)} Assume in addition $\delta >0,$ that is, $\theta <\gamma
\delta <0.$ Then (\ref{pg}), (\ref{coef}), (\ref{delta}) yield $p>0,$ $%
\Delta <0,$ $a>0$ and $q>0.$ The equilibrium point $E_{3}$ is an unstable
node in this case, whenever it exists in $Q1.$ The bifurcation diagram is
depicted in Figure \ref{d1} (I).\vspace{6pt}

\begin{figure}[h!]
\centering
\includegraphics[width=0.9\textwidth]{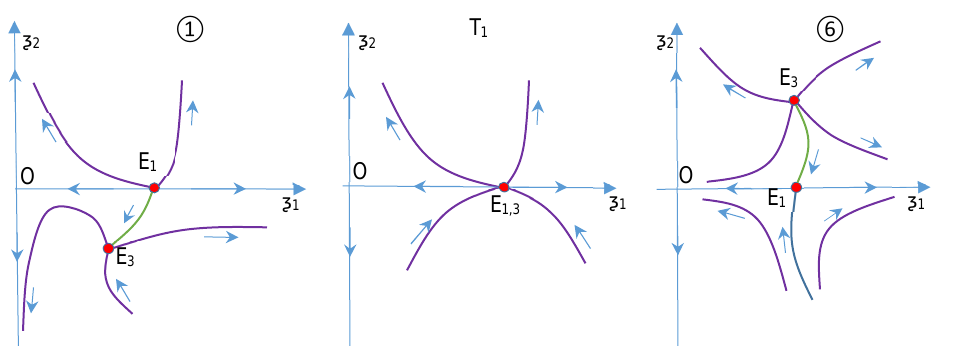}
\caption{The behavior of the system (\protect\ref{s2}) when $\left( \protect%
\mu _{1},\protect\mu _{2}\right)$ crosses $T_1$ in the first bifurcation diagram (I). }
\label{midT1}
\end{figure}

Let us describe the phase portraits on the bifurcation curves $T_{1,2},$ $%
X_{+,-}$ and $Y_{+,-}$ in the first diagram (I). We will exemplify this on
the curve $T_{1}.$ The description is similar for the other curves and the
next diagrams. Since $T_{1}$ is a transcritical bifurcation curve and the
colliding points $E_{3}$ and $E_{1}\left( \frac{1}{\theta }\mu
_{1}+\allowbreak O\left( \mu _{1}^{2}\right) ,0\right) $ have the
eigenvalues $0$ and $-\mu _{1}+O\left( \mu _{1}^{2}\right) ,$ the phase
portrait of the system (\ref{s2}) is known when $\left( \mu _{1},\mu
_{2}\right) \in T_{1}$ and given in Figure \ref{midT1}(T1). Notice that $%
\theta <0$ and, from $\lambda _{1}\lambda _{2}=-\left( \theta -\gamma \delta
\right) \xi _{1}\xi _{2}\left( 1+O\left( \left\vert \mu \right\vert \right)
\right) $ given in (\ref{l12}), $E_{3}\left( \xi _{1},\xi _{2}\right) $ is a
saddle when $\left( \mu _{1},\mu _{2}\right) $ lies in the region 1 from
Figure \ref{d1} (I), since $\xi _{1}>0$ and $\xi _{2}<0,$ see Figure \ref%
{midT1}(1). After $\left( \mu _{1},\mu _{2}\right) $ crosses $T_{1}$ into
the region 6, $E_{3}\left( \xi _{1},\xi _{2}\right) $ becomes an unstable
node and lies in the first quadrant, Figure \ref{midT1}(6).

Comparing the two phase portraits from Figures \ref{midT1}(1)-\ref{midT1}%
(T1), we observe they coincide in the first quadrant. As a conclusion we
have the next remark.

\begin{remark}
Denote by $T_{1}^{vir}$ the region to the left or right of $T_{1}$ where $%
E_{3}$ is virtual (i.e. $E_{3}\notin Q1).$ For example, $T_{1}^{vir}$ is the
region 1 in Figure \ref{d1} (I). Then, the behaviour of the
system (\ref{s2}) when $\left( \mu _{1},\mu _{2}\right) \in T_{1}$ coincides
in the first quadrant $Q1$ to the one for $\left( \mu _{1},\mu _{2}\right)
\in T_{1}^{vir}.$ In other words, the phase portrait on $T_{1}$ is the same
(in Q1) as the one to the left or right of $T_{1}$ where $E_{3}$ is virtual.
The same scenario applies for $T_{2}$ and all the other bifurcation curves $%
X_{+,-}$ and $Y_{+,-},$ with the difference that $E_{3}$ is replaced
correspondingly with $E_{1}$ or $E_{2},$ since in these latter cases, $E_{1}$
or $E_{2}$ collide to $O.$ This rule applies for all bifurcation diagrams
from this section.
\end{remark}

\begin{remark}
On the curves $C_{1}-C_{4},$ which separate a node from a focus, $E_{3}$ is
a node (since $q=0$), thus, the phase portrait on each of these curves
coincides to the phase portrait on the left or right of the curves where $%
E_{3}$ is a node. 
\end{remark}

\begin{figure}[h!]
\centering
\includegraphics[width=0.45\textwidth]{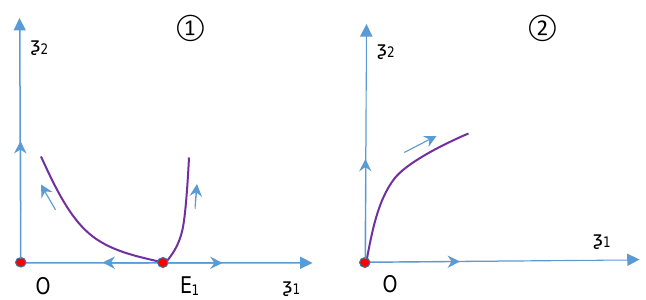} %
\includegraphics[width=0.45\textwidth]{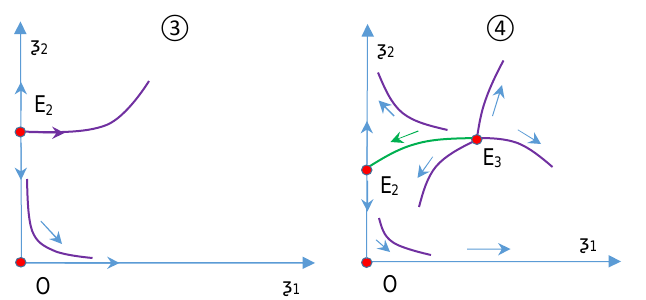} %
\includegraphics[width=0.45\textwidth]{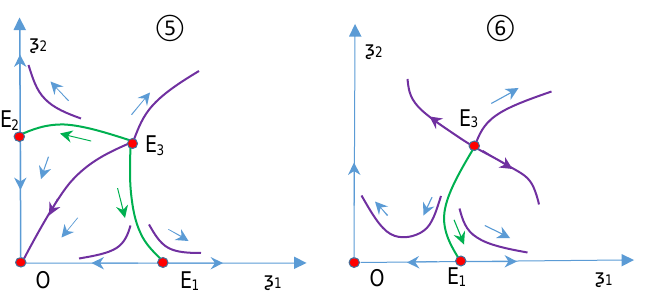} %
\includegraphics[width=0.45\textwidth]{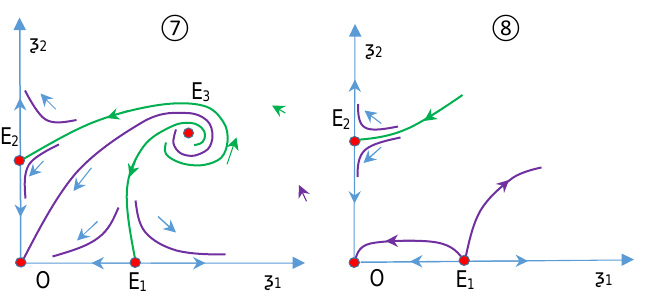} %
\includegraphics[width=0.45\textwidth]{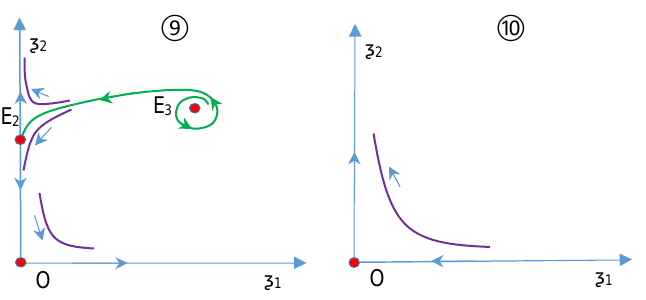} %
\includegraphics[width=0.45\textwidth]{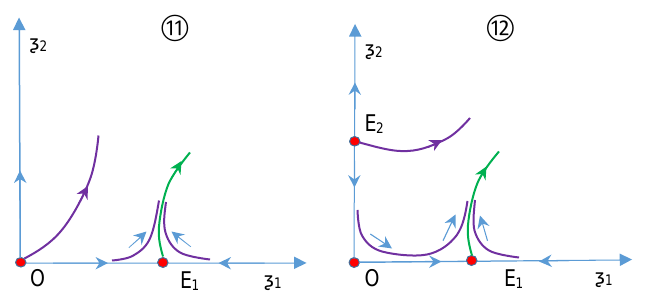} %
\includegraphics[width=0.45\textwidth]{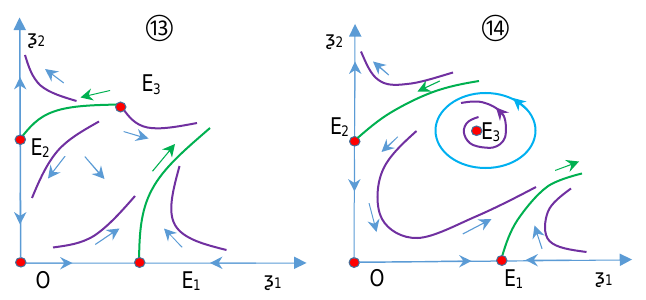} %
\includegraphics[width=0.45\textwidth]{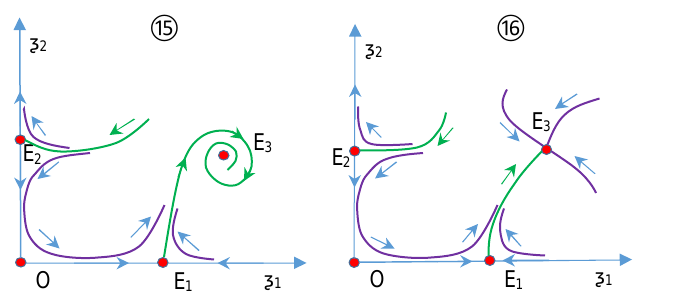} %
\includegraphics[width=0.45\textwidth]{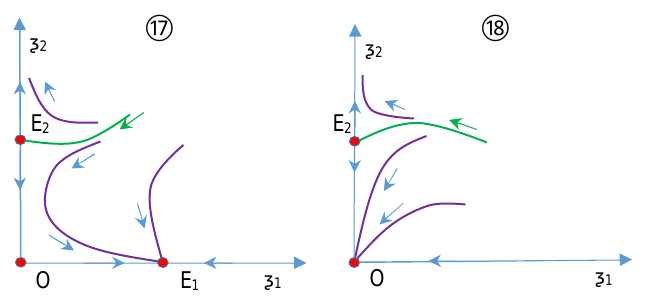} %
\includegraphics[width=0.45\textwidth]{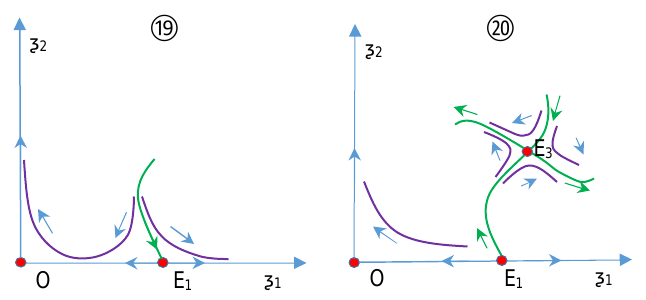}
\caption{Phase portraits corresponding to different regions of the
bifurcation diagrams.}
\label{pp1}
\end{figure}

\begin{figure}[h!]
\centering
\includegraphics[width=0.45\textwidth]{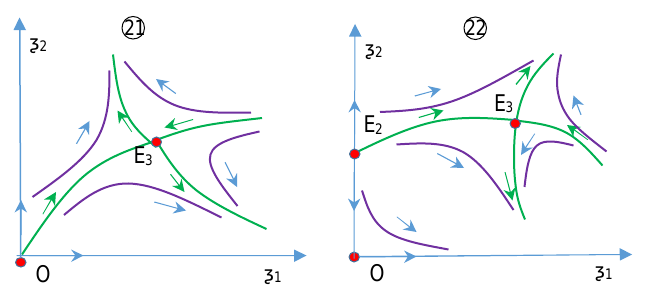} %
\includegraphics[width=0.45\textwidth]{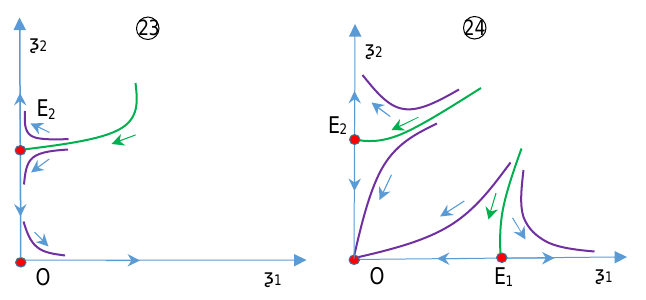} %
\includegraphics[width=0.45\textwidth]{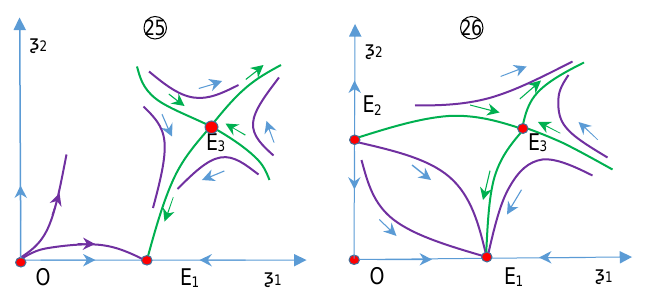}
\caption{Phase portraits corresponding to different regions of the
bifurcation diagrams.}
\label{pp2}
\end{figure}

\textbf{Case 1.b)} Consider $\theta -\gamma \delta <0,$ $\delta <0,$ and $%
\theta <0.$ Therefore $p>0$ which, in turn, implies that a Hopf bifurcation
cannot occur at $E_{3}.$ From $\Delta >0$ and Lemma \ref{lemma2}, two
bifurcation curves arise from $q=0,$ namely 
\begin{equation}
C_{1}:\mu _{2}=e_{1}\mu _{1}\text{ and }C_{2}:\mu _{2}=e_{2}\mu _{1},
\label{c12}
\end{equation}%
where $a,b,c$ are given by (\ref{coef}). The two curves $C_{1,2}$ control
the sign of $q.$ Denote by $m_{C_{1}}=e_{1},$ $m_{C_{2}}=e_{2},$ $m_{T_{1}}=%
\frac{\delta }{\theta }$ and $m_{T_{2}}=\frac{1}{\gamma }$ the slopes of the
straight lines $C_{1},$ $C_{2},$ $T_{1}$ (\ref{T1}) and $T_{2}$ (\ref{T2})
respectively. Notice that $m_{T_{1}}>0$ and $m_{T_{2}}<0.$ To determine the
relative position of these slopes on the real axis, we need to know the
signs of the following expressions: 
\begin{equation}
P_{1}=m_{1}\cdot m_{2}=\frac{\left( \theta -\gamma \delta \right) ^{2}}{%
\gamma ^{2}a},~S_{1}=m_{1}+m_{2}=\frac{-4\left( \theta -\gamma \delta
\right) \theta }{a}\left( \frac{\gamma +1}{2\gamma }-\gamma \frac{\delta }{%
\theta }\right)  \label{ps1}
\end{equation}%
and%
\begin{equation}
P_{2}=m_{1}^{\prime }\cdot m_{2}^{\prime }=\frac{\left( \theta -\gamma
\delta \right) ^{2}}{a},~S_{2}=m_{1}^{\prime }+m_{2}^{\prime }=\frac{4\left(
\theta -\gamma \delta \right) \theta \gamma }{a}\left( \frac{\gamma +1}{%
2\gamma }-\frac{\delta }{\theta }\right) ,  \label{ps2}
\end{equation}%
where $m_{i}=e_{i}-\frac{1}{\gamma }$ and $m_{i}^{\prime }=e_{i}-\frac{%
\delta }{\theta },$ $i=1,2.$

Taking into account that the sign of $a=4\gamma ^{2}\theta ^{2}\left[ \left( 
\frac{\gamma +1}{2\gamma }\right) ^{2}-\frac{\delta }{\theta }\right] $ is
undetermined in this subcase, we obtain the following possible situations.$~$%
\newline
$(i)$ Consider $a>0,$ that is, $0<\frac{\delta }{\theta }<\left( \frac{%
\gamma +1}{2\gamma }\right) ^{2}$ and $\gamma \neq -1.$

If $-1<\gamma <0,$ then from (\ref{ps2}) we deduce $P_{2}>0$ and $S_{2}>0.$
Therefore we have 
\begin{equation*}
\frac{1}{\gamma }<0<\frac{\delta }{\theta }<e_{1}<e_{2}.
\end{equation*}%
The bifurcation diagram is depicted in Figure \ref{d1} (II).

\begin{figure}[h!]
\centering
\includegraphics[width=0.9\textwidth]{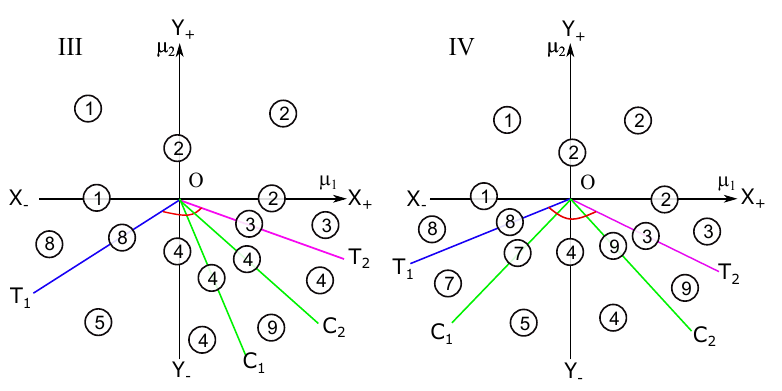}
\caption{Bifurcation diagrams of system (\protect\ref{s2}) for: $\protect%
\theta-\protect\gamma\protect\delta<0,$ $\protect\delta<0,$ $\protect\theta%
<0,$ and $a>0$ and $\protect\gamma<-1$ (III), respectively, $\protect\theta-%
\protect\gamma\protect\delta<0,$ $\protect\delta<0,$ $\protect\theta<0$ and $%
a<0$ (IV).}
\label{d2}
\end{figure}

If $\gamma <-1,$ then from (\ref{ps1}) we deduce $P_{1}>0$ and $S_{1}<0.$
Therefore we have 
\begin{equation*}
e_{1}<e_{2}<\frac{1}{\gamma }<0<\frac{\delta }{\theta }.
\end{equation*}%
The bifurcation diagram is depicted in Figure \ref{d2} (III).

$(ii)$ Assume $a<0,$ that is, $\frac{\delta }{\theta }>\left( \frac{\gamma +1%
}{2\gamma }\right) ^{2}$ for all $\gamma <0.$ Hence $e_{2}<0<e_{1}.$
Moreover, $P_{1}<0$ and $P_{2}<0,$ whence 
\begin{equation*}
e_{2}<\frac{1}{\gamma }<0<\frac{\delta }{\theta }<e_{1}.
\end{equation*}%
The bifurcation diagram is depicted in Figure \ref{d2} (IV).\newline
$(iii)$ Finally, let $a=0,$ that is, $\theta =\frac{4\gamma ^{2}\delta }{%
(1+\gamma )^{2}};$ $\gamma \neq -1$ in this case because, otherwise, $%
a=-4\theta \delta \neq 0$ at $\gamma =-1.$ Hence $b=\frac{2\delta ^{2}\gamma
^{2}(1-\gamma )^{3}}{(1+\gamma )^{3}}\not=0,$ while $c\neq 0$ for all $%
\theta \delta \neq 0$ and $\gamma <0.$ The sign of $q$ cannot be concluded
in this case from (\ref{q}) because a term in $\mu _{2}^{n}$ with $n\geq 2$
is missing. Therefore, we have to use terms up to order three in the general
expression of $q$ given by (\ref{qg}) in order to determine its leading
(lowest) terms. At $\theta =\frac{4\gamma ^{2}\delta }{(1+\gamma )^{2}},$ we
obtain from (\ref{qg}) that $q$ in its lowest terms reads 
\begin{equation}
q=\allowbreak c\mu _{1}^{2}\left( 1+O_{1}\left( \left\vert \mu \right\vert
\right) \right) -2b\mu _{1}\mu _{2}\left( 1+O_{2}\left( \left\vert \mu
\right\vert \right) \right) +d\mu _{2}^{3}\left( 1+O_{3}\left( \left\vert
\mu \right\vert \right) \right) ,  \label{qa0}
\end{equation}%
where $c=\frac{4\gamma ^{2}+3\gamma +1}{4\gamma ^{2}\left( 1-\gamma \right) }%
,$ $b=-\frac{\gamma +1}{2\left( \gamma -1\right) },$ $d=\frac{4\gamma \left(
\gamma +1\right) }{\delta ^{2}\left( \gamma -1\right) ^{6}}\left( \left(
2S-N\right) \gamma ^{4}+d_{31}\gamma ^{3}+d_{21}\allowbreak \gamma
^{2}+d_{11}\gamma +M\delta -N\right) ,$ respectively, $d_{31}=6S-4N+3M\delta
,$ $d_{21}=16P\delta ^{2}+7M\delta -6N+6S$ and $d_{11}=2S-4N+5M\delta .$

The sign of $q$ from (\ref{qa0}) can be determined. Applying a method based
on IFT as in Lemma \ref{lemma2}, the equation $q=0$ is satisfied along two
bifurcation curves, given by%
\begin{equation}
C_{3}=\left\{ \left( \mu _{1},\mu _{2}\right) \left\vert \mu _{1}=\frac{2b}{c%
}\mu _{2}+O\left( \mu _{2}^{2}\right) \right. \right\}  \label{c3}
\end{equation}%
and 
\begin{equation}
C_{4}=\left\{ \left( \mu _{1},\mu _{2}\right) \left\vert \mu _{1}=-d\frac{%
\gamma -1}{\gamma +1}\mu _{2}^{2}+O\left( \mu _{2}^{3}\right) \right.
\right\}  \label{c4}
\end{equation}%
for all $\gamma \neq -1.$ The two curves determine in the parametric plane
the sign of $q$ for all $\left\vert \mu \right\vert $ sufficiently small.
Notice that $C_{4}$ is tangent to the $\mu _{2}-$axis at $O.$

\begin{figure}[h!]
\centering
\includegraphics[width=0.9\textwidth]{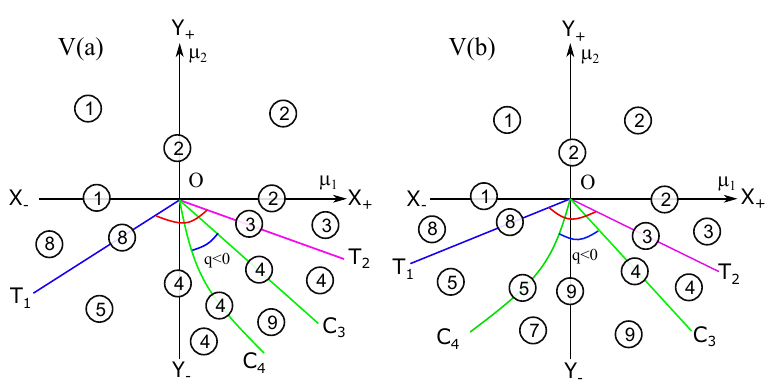}
\caption{Bifurcation diagrams of system (\protect\ref{s2}) when $\protect%
\theta-\protect\gamma\protect\delta <0,$ $\protect\delta <0,$ $\protect\theta%
<0,$ $a=0,$ $\protect\gamma<-1,$ respectively, $d<0$ on V(a) and $d>0$ on
V(b).}
\label{d3a}
\end{figure}

\begin{figure}[h!]
\centering
\includegraphics[width=0.9\textwidth]{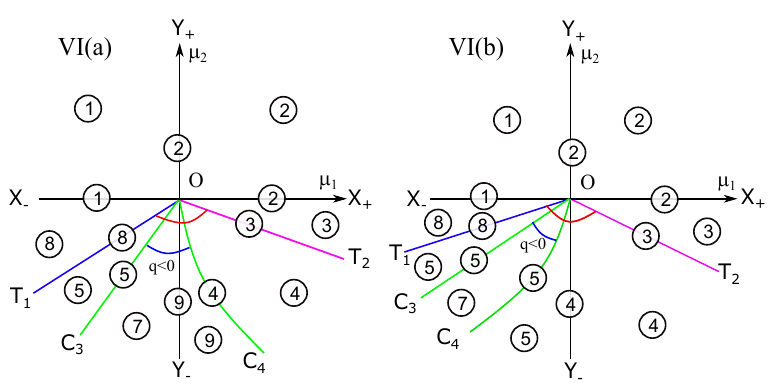}
\caption{Bifurcation diagrams of system (\protect\ref{s2}) when $\protect%
\theta-\protect\gamma\protect\delta <0,$ $\protect\delta <0,$ $\protect\theta%
<0,$ $a=0,$ $-1<\protect\gamma<0,$ respectively, $d>0$ on VI(a) and $d<0$ on
VI(b). }
\label{d3b}
\end{figure}

Two cases arise, namely $b<0$ (iff $\gamma <-1$) and then $\frac{1}{\gamma }>%
\frac{c}{2b}$ (see Figures \ref{d3a}, V(a)-(b)) and $b>0$ (iff $-1<\gamma <0$%
) and then $\frac{\delta }{\theta }<\frac{c}{2b}$ (see Figures \ref{d3b},
VI(a)-(b)), respectively. Notice that $\frac{c}{2b}=\allowbreak \frac{%
4\gamma ^{2}+3\gamma +1}{4\left( \gamma +1\right) \gamma ^{2}},$
respectively, $\frac{1}{\gamma }-\frac{c}{2b}=\allowbreak \frac{\gamma -1}{%
4\gamma ^{2}\left( \gamma +1\right) }$ and $\frac{\delta }{\theta }-\frac{c}{%
2b}=\allowbreak \frac{\gamma -1}{4\left( \gamma +1\right) }.$ In each of the
two cases we considered that $d$ can be either positive or negative. When $%
d=0,$ a further analysis is required to determine the sign of $q.$ This will
need all terms in the general expression of $q$ from (\ref{qg}).

\begin{figure}[h!]
\centering
\includegraphics[width=0.9\textwidth]{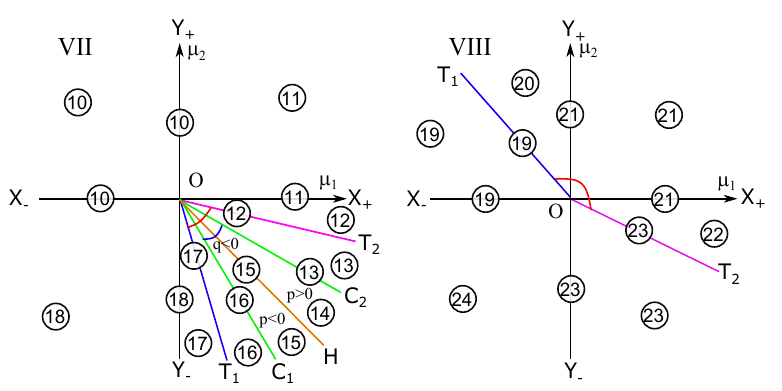}
\caption{Parametric portrait of system (\protect\ref{s2}) when $\protect%
\theta - \protect\gamma \protect\delta <0,$ $\protect\delta <0,$ and $%
\protect\theta>0 $ (VII), respectively, $\protect\theta - \protect\gamma 
\protect\delta >0, $ $\protect\delta >0,$ and $\protect\theta<0$ (VIII). }
\label{d4}
\end{figure}

\begin{figure}[h!]
\centering
\includegraphics[width=0.9\textwidth]{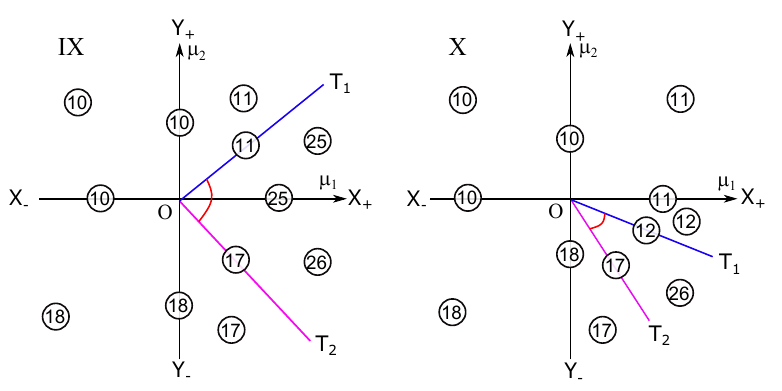}
\caption{Parametric portrait of system (\protect\ref{s2}) when $\protect%
\theta - \protect\gamma \protect\delta >0,$ $\protect\delta >0,$ and $%
\protect\theta>0 $ (IX), respectively, $\protect\theta - \protect\gamma 
\protect\delta >0$ and $\protect\delta <0$ (X).}
\label{d5}
\end{figure}

\bigskip

\textbf{Case 1.c)} Consider $\theta -\gamma \delta <0,$ $\delta <0,$ and $%
\theta >0.$ Thus $0<\theta <\gamma \delta .$ In this case, a Hopf
bifurcation occurs at $E_{3}.$ Indeed, $p=0$ along the bifurcation curve 
\begin{equation*}
H=\left\{ \left( \mu _{1},\mu _{2}\right) \left\vert \mu _{2}=\allowbreak
\allowbreak \frac{\theta -\delta }{\theta \left( \gamma -1\right) }\mu
_{1}+O\left( \mu _{1}^{2}\right) ,\,\theta \neq \delta \right. \right\} \,
\end{equation*}%
and $\left. \frac{\partial p}{\partial \mu _{2}}\right\vert _{H}=\allowbreak 
\frac{\left( \gamma -1\right) \theta }{2\left( \theta -\gamma \delta \right) 
}\neq 0.$ Notice that $\theta \neq \delta $ is satisfied in this case. Then $%
q=\mu _{1}^{2}\frac{\theta -\gamma \delta }{\theta \left( \gamma -1\right)
^{2}}<0$ on $H$ and the eigenvalues associated to the equilibrium point $%
E_{3}$ are $\lambda _{1,2}=\pm i\sqrt{-q}$ along the curve $H.$ The
following order is obtained in this case, namely: 
\begin{equation}
\frac{\delta }{\theta }<e_{1}<\frac{\theta -\delta }{\theta \left( \gamma
-1\right) }<e_{2}<\frac{1}{\gamma }<0\,.  \label{oh}
\end{equation}

Indeed, from $\gamma <0$ we have that $\Delta =-4\gamma \delta \left( \theta
-\gamma \delta \right) ^{3}>0,$ $a>0$ and the equation $a\mu _{2}^{2}-2b\mu
_{1}\allowbreak \mu _{2}+c\allowbreak \mu _{1}^{2}=0$ in $\mu _{2}$ has two
roots $\mu _{2}=e_{1}\mu _{1}$ and $\mu _{2}=e_{2}\mu _{1},$ where $e_{1}=%
\frac{b-\sqrt{\Delta }}{a}$ and $e_{2}=\frac{b+\sqrt{\Delta }}{a},$ with $%
e_{1}<e_{2}.$ Therefore we deduce $m_{1}\cdot m_{2}>0,$ $m_{1}^{\prime
}\cdot m_{2}^{\prime }>0$ and $\left( e_{1}-\frac{\theta -\delta }{\theta
\left( \gamma -1\right) }\right) \left( e_{2}-\frac{\theta -\delta }{\theta
\left( \gamma -1\right) }\right) =\frac{4\left( \theta -\gamma \delta
\right) ^{3}}{\theta (\gamma -1)^{2}a}<0.$ Also we have 
\begin{equation*}
\theta +\theta \gamma -2\gamma \delta <\theta \gamma -\theta <0\text{ and }%
\theta +\theta \gamma -2\gamma ^{2}\delta >\theta -\theta \gamma >0.
\end{equation*}%
Hence $m_{1}^{\prime }+m_{2}^{\prime }=\frac{\allowbreak 2}{a}\left( \theta
-\gamma \delta \right) \left( \theta +\theta \gamma -2\gamma \delta \right)
>0,$ $m_{1}+m_{2}=\allowbreak -2\frac{\theta -\gamma \delta }{\gamma a}%
\left( \theta +\theta \gamma -2\gamma ^{2}\delta \right) <0,$ and the
conclusion follows. The bifurcation diagram is depicted in Figure \ref{d4}
(VII). When $l_{1}\left( 0\right) <0,$ a limit cycle arises for $\mu
_{2}>\allowbreak \allowbreak \frac{\theta -\delta }{\theta \left( \gamma
-1\right) }\mu _{1}.$

If the slope of $H$ is denoted by $m_{H}=\frac{\theta -\delta }{\theta
\left( \gamma -1\right) }$ and $m_{C_{1}}=e_{1},$ $m_{C_{2}}=e_{2},$ $%
m_{T_{1}}=\frac{\delta }{\theta }$ and $m_{T_{2}}=\frac{1}{\gamma }$
represent the slope of the straight line $C_{1},$ $C_{2},$ $T_{1},$ $T_{2},$
respectively, then (\ref{oh}) leads to the following order for the slopes of
the curves: 
\begin{equation*}
m_{T_{1}}<m_{C_{1}}<m_{H}<m_{C_{2}}<m_{T_{2}}<0\,.
\end{equation*}

\bigskip

\textbf{Case 2. }In the following we consider $\theta -\gamma \delta >0.$
Hence, from (\ref{l12}), we have $\lambda _{1}\lambda _{2}<0.$ Therefore the
equilibrium point $E_{3}$ is a saddle, provided that it exists.\vspace{6pt}

\textbf{Case 2.a)} Consider $\theta -\gamma \delta >0,$ $\delta >0,$ and $%
\theta <0.$ Then $p>0$ and there are no Hopf bifurcation. But we have $%
\Delta >0,$ $a>0.$ The equilibrium point $E_{3}$ exists if and only if $%
\delta \mu _{1}-\theta \mu _{2}>0$ and $\mu _{1}-\gamma \mu _{2}>0.$ The
bifurcation diagram is depicted in Figure \ref{d4} (VIII).\vspace{6pt}

\textbf{Case 2.b)} Consider $\theta -\gamma \delta >0,$ $\delta >0,$ and $%
\theta >0.$ Then $\gamma \delta <0<\theta $ and $\Delta >0.$ The bifurcation
diagram is depicted in Figure \ref{d5} (IX).\vspace{6pt}

\textbf{Case 2.c)} Consider $\theta -\gamma \delta >0$ and $\delta <0.$ Then 
$\theta >\gamma \delta >0,$ which leads to $\Delta <0$ and $a>0.$ The
bifurcation diagram is depicted in Figure \ref{d5} (X).\vspace{6pt}

In conclusion, we have obtained \textit{twelve bifurcation diagrams}
corresponding to all possible situations. The Hopf bifurcation and,
consequently, the existence of limit cycles emerging from the bifurcation
are only possible in the bifurcation diagram VII (Figure \ref{d4}),
corresponding to the Case 1.c). All possible types of the four equilibrium
points arising in the diagrams are summarized in Table \ref{table:1} and
Table \ref{table:2}. There are \textit{twenty six} distinct generic phase
portraits appearing in the bifurcation diagrams.

\begin{table}[h!]
\centering
\begin{tabular}{c|ccccccccccccc}
& $1$ & $2$ & $3$ & $4$ & $5$ & $6$ & $7$ & $8$ & $9$ & $10$ & $11$ & $12$ & 
$13$ \\ \hline
$O$ & s & un & s & s & sn & s & sn & sn & s & s & un & s & s \\ 
$E_{1}$ & un & - & - & - & s & s & s & un & - & - & s & s & s \\ 
$E_{2}$ & - & - & un & s & s & - & s & s & s & - & - & un & s \\ 
$E_{3}$ & - & - & - & un & un & un & uf & - & uf & - & - & - & un%
\end{tabular}%
\caption{\textit{The types of the equilibrium points of system (\protect\ref%
{s2}) on different regions of the bifurcation diagrams; the abbreviations s,
sn, un, sf, and uf stand for saddle, stable node, unstable node, stable
focus, and unstable focus, respectively.}}
\label{table:1}
\end{table}

\begin{table}[h]
\centering
\begin{tabular}{c|ccccccccccccc}
& $14$ & $15$ & $16$ & $17$ & $18$ & $19$ & $20$ & $21$ & $22$ & $23$ & $24$
& $25$ & $26$ \\ \hline
$O$ & s & s & s & s & sn & s & s & un & s & s & sn & un & s \\ 
$E_{1}$ & s & s & s & sn & - & s & un & - & - & - & s & sn & sn \\ 
$E_{2}$ & s & s & s & s & s & - & - & - & un & s & s & - & un \\ 
$E_{3}$ & uf & sf & sn & - & - & - & s & s & s & - & - & s & s%
\end{tabular}%
\caption{\textit{Continuation of Table \protect\ref{table:1}.}}
\label{table:2}
\end{table}

\section{Conclusions}

A conclusion can be drawn from Theorems \ref{t1}-\ref{t2} and the above
twelve bifurcation diagrams. It is given in the next theorem. Its proof is a
direct consequence of the fact that, in all the above results from this
section (with the exception of the case when a Hopf bifurcation occurs), we
only needed terms up to order two in (\ref{s2}). These include also the case 
$a=0$ because the curves $C_{3,4}$ do not change qualitatively the behavior
of the system (the node-focus topological equivalence).

\begin{theorem}
Assume $\theta \delta \neq 0$ and $\theta -\gamma \delta \neq 0.$ If $\theta
<0$ or $0<\gamma \delta <\theta $ or $\gamma \delta <0<\theta ,$ the system (%
\ref{s2}) is locally topologically equivalent near the origin to 
\begin{equation*}
\left\{ 
\begin{tabular}{lll}
$\frac{d\xi _{1}}{dt}$ & $=$ & $\xi _{1}\left( \mu _{1}-\theta \xi
_{1}+\gamma \xi _{2}\right) $ \\ 
$\frac{d\xi _{2}}{dt}$ & $=$ & $\xi _{2}\left( \mu _{2}-\delta \xi _{1}+\xi
_{2}\right) $%
\end{tabular}%
\right. .
\end{equation*}
\end{theorem}

\section{Acknowledgments}

This research was supported by Horizon2020-2017-RISE-777911 project. We
thank to Prof. Jaume Llibre for his useful suggestions in elaborating the
article and to the anonymous referee for his valuable comments.

\end{document}